\documentclass[10pt,reqno]{amsart}

\usepackage{amssymb}

\usepackage{amsmath}

\usepackage{bm}

\usepackage{verbatim}

\makeatletter

\@addtoreset{equation}{section}    

\makeatother

\newtheorem{thm}{Theorem}[section]

\newcommand{\T}{\mathbb T}

\newlength{\defbaselineskip}

\newcommand{\setlinespacing}[1]%

\title[K${\bf d}$V Preserves White Noise]
{K${\bf d}$V Preserves White Noise}

\author{Jeremy Quastel}
\author{Benedek Valk\' o}

\address{\noindent Departments of Mathematics and Statistics, University of
Toronto
\newline e-mail:  \rm \texttt{quastel@math.toronto.edu, valko@math.toronto.edu }}

\date{\today}

\begin{document}

\maketitle

\begin{abstract}
It is shown that white noise is an invariant measure for the Korteweg-deVries equation
on $\mathbb T$.  This is a consequence of recent results of Kappeler and Topalov establishing the well-posedness of the equation on appropriate negative Sobolev spaces, together with a result of
  Cambronero and McKean that white noise is the image under the Miura transform (Ricatti map) of the (weighted) Gibbs measure for
the modified KdV equation, proven to be invariant
for that equation by
Bourgain.\end{abstract}

\section{K${\rm d}$V on $H^{-1}(\Bbb T)$ and White Noise}

The Korteweg-deVries equation (KdV) on $\mathbb T= \mathbb R/\mathbb Z$,
\begin{equation}\label{kdv}
u_t - 6 u u_x +u_{xxx} =0, \qquad u(0) = f
\end{equation}
defines  nonlinear evolution operators
\begin{equation}
{\mathcal S}_t f = u(t)
\end{equation}
$-\infty<t < \infty$ on smooth functions $f:\mathbb T\to \mathbb R$.

\begin{thm} \label{kt1} (Kappeler and Topalov \cite{KT1}) $\mathcal S_t$ extends to a continuous group of nonlinear evolution operators
\begin{equation}
\bar{\mathcal S}_t: H^{-1}(\mathbb T) \to H^{-1}(\mathbb T).
\end{equation} \end{thm} In concrete terms, take $f\in H^{-1}(\mathbb T)$ and let $f_N$ be smooth functions on 
$\mathbb T$ with $\|f_N- f\|_{H_{-1}(\mathbb T)}\to 0$ as
$N\to \infty$.  
Let $u_N(t)$ be the (smooth) solutions of (\ref{kdv}) with initial data $f_N$.  Then there is a unique
$u(t)\in H^{-1}(\mathbb T)$ which we call
$u(t)= \bar{\mathcal S}_t f$ with $\|u_N(t) - u(t)\|_{H^{-1}(\mathbb T)}\to 0$.

White noise on $\Bbb T$ is the unique probability
measure $Q$ on the space $\mathcal D(\Bbb T)$ of
distributions on $\Bbb T$ satisfying
\begin{equation}
\int e^{i \langle\lambda, u\rangle} dQ(u) = e^{-\frac12 \|\lambda\|_2^2}
\end{equation}
for any smooth function $\lambda$ on $\mathbb T$
where $\|\cdot\|_2^2 =\langle\cdot, \cdot\rangle$
are the $L^2(\Bbb T, dx)$ norm and inner product
(see \cite{H}).  

Let $\{e_n\}_{n=0,1,2,\ldots}$ be an orthonormal
basis of smooth functions in $L^2(\mathbb T)$ with $e_0=1$.  White noise
 is represented as  $u=\sum_{n=0}^\infty x_n e_n$ where and $x_n$ are independent Gaussian random variables, each with mean $0$ and variance $1$.
Hence $Q$ is supported in $H^{-\alpha}(\Bbb T)$ for any
$\alpha>1/2$.    

Mean zero white noise $Q_0$ on $\mathbb T$ is the probability measure on 
distributions $u$ with $\int_{\mathbb T} u=0$ satisfying
\begin{equation}
\int e^{i \langle\lambda, u\rangle} dQ_0(u) = e^{-\frac12 \|\lambda\|_2^2}
\end{equation}
for any mean zero smooth function $\lambda$ on $\mathbb T$.  It is represented as
$u=\sum_{n=1}^\infty x_ne_n$.

Recall that if $f:X_1\to X_2$ is a measurable map
between metric spaces and $Q$ is a probability measure on $(X_1,\mathcal B(X_1))$, then the 
{\it pushforward} $f^*Q$ is the measure on 
$X_2$ given by $f^*Q (A)= Q( \{ x~:~
f(x) \in A\})$ for any Borel set $A\in \mathcal B (X_2)$.

Our main result is:

\begin{thm}\label{main}  White noise $Q_0$ is invariant under KdV;  for any $t\in \Bbb R$,
\begin{equation}\bar{\mathcal S}_t^{*} Q_0 = Q_0.\end{equation}  \end{thm}

{\it Remarks.}  1. In terms of classical solutions
of KdV, the meaning of Theorem \ref{main} is as follows.  Let $f_N$, $N=1,2,\ldots$ be a sequence
of smooth mean zero random initial data approximating mean zero white 
noise.  For example, one could take $f_N(\omega) =\sum_{n=1}^N x_n(\omega) e_n$ where $x_n$ and $e_n$ are as above.  Solve the
KdV equation for each $\omega$ up to a fixed time $t$ to obtain
$\mathcal S_t f_N$.  The limit in $N$ exists \cite{KT1} in $H^{-1}(\Bbb T)$, for almost every value of $\omega$ , and is again a 
white noise.

2.  It follows immediately that 
\begin{equation}
\hat{\mathcal S}_t: L^2(Q_0)\to L^2(Q_0),  \quad
(\hat{\mathcal S}_t \Phi) (f) = \Phi(\bar{\mathcal S}_t f)
\end{equation}
are a group of {\it unitary} transformations of $L^2(Q_0)$, defining a continuous Markov process $u(t)$, $t\in (-\infty, \infty)$ on $H^{-1}(\Bbb T)$ with 
Gaussian white noise one dimensional marginals, invariant under time$+$space inversions.  The correlation functions $S(x,t) = 
\int f(0)\bar{\mathcal S}_t f (x) 
dQ_0$ may  have
an interesting structure.

3.  A similar result holds without the mean zero
condition, but now the mean $m= \int_{\mathbb T} u$ is distributed not as
an independent Gaussian, but as one conditioned
to have $m\ge -\lambda_0(u)$ where $\lambda_0(u)$
is the principal eigenvalue of $-\frac{d^2}{dx^2}
+ u$.  Since the  addition of constants produces a trivial
rotation in the KdV equation it seems more natural
to consider the mean zero case.

4.  $Q_0$ is certainly {\it not} the only invariant measure for 
KdV.  The Gibbs measure formally written as
$Z^{-1}1\left(\int_{\Bbb T}u^2\le K\right)e^{-{\mathcal H}_2}$ where 
\begin{equation} \label{h2}
{\mathcal H}_2(u)= -\int_{\Bbb T} u^3 -\frac12
u_x^2
\end{equation} is known to be invariant \cite{Bo}.   Note that (after subtraction of the mean) this Gibbs measure is supported on a set of $Q_0$-measure $0$.  
$Q_0$ is also a Gibbs measure, corresponding to
the Hamiltonian \begin{equation}
\label{h1}
{\mathcal H}_1(u)=\int_{\Bbb T} u^2.
\end{equation}  The existence of two Gibbs measures corresponds to the
bihamiltonian structure of KdV:  It can be written
\begin{equation}
\dot u = J_i \frac{\delta {\mathcal H}_i}{\delta u}, \qquad i=1,2
\end{equation}
with 
symplectic forms $J_1= \partial_x^3 
+4u\partial_x +2\partial_x u$ and
$J_2= \partial_x$.
Because of all the conservation laws of KdV, there are  many other invariant measures as well.

5.  We were led to Theorem \ref{main} after noticing that the discretization of KdV used by
Kruskal and Zabusky in the numerical investigation of solitons,
\begin{equation}\label{sp}
\dot u_i= (u_{i+1}+u_i+u_{i-1}) (u_{i+1}-u_{i-1})
-(u_{i+2}-2u_{i+1}+2u_{i-1}-u_{i-2}),
\end{equation}
preserves discrete white noise (independent Gaussians mean $0$ and variance $\sigma^2>0$).   The invariance follows from two simple
properties of the special discretization (\ref{sp}).  First of all  $\dot u_i = b_i$
preserves Lebesgue measure whenever $\nabla\cdot b
=\sum_i \partial_ib_i = 0$, and  (\ref{sp}) is of this
form.  Furthermore, it is easy to check (though something of a miracle) that $\sum_i u_i^2$ is
invariant under (\ref{sp}).  Hence $Z^{-1} e^{-\frac{1}{2\sigma^2} \sum_i u_i^2}\prod du_i$ is also invariant.

 Note that the discretization (\ref{sp}) is {\it not}
completely integrable, and we are not aware of a completely 
integrable discretization which does conserve discrete
white noise.  For example, consider the following  family of completely integrable discretizations of KdV, depending on 
a real parameter $\alpha$ \cite{al};
\begin{eqnarray}
 \dot u_i & = &  (1-\alpha u_i) \{-\alpha u_{i-1}(u_{i-2}-u_{i})
 -\alpha(u_{i-1}+ 2u_i+u_{i+1})(u_{i-1}-u_{i+1})
 \nonumber\\ &&
 -\alpha u_{i+1}(u_{i}-u_{i+2}) + u_{i-2}-2u_{i-1}+2u_{i+1}-u_{i+2}\}.
 \end{eqnarray}
 They conserve Lebesgue measure by the Liouville theorem.  We want $\alpha\neq 0$; otherwise the quadratic term of KdV
 is not represented.  In that case the conserved quantity
 analagous to $\int_{\Bbb T}u^2$ is
 \begin{equation}{\mathcal Q}=\sum_i u_i^2+2 u_i u_{i+1}.
 \end{equation} 
But $\mathcal Q$ is non-definite, and hence the corresponding
 measure $e^{-\mathcal Q}\prod_i du_i$ cannot be normalized 
 to make a probability measure.

6.  At a completely formal level the proof proceeds as follows.  Note first of all that the flow generated by $u_t= u_{xxx}$ is easily solved and 
seen to preserve white noise.  So  consider the
Burgers' flow $u_t=2u u_x$.
\begin{equation}
\partial_t \int f(u(t))e^{-\int u^2}
 =  \int \langle \frac{\delta f}{\delta u},  u_t\rangle 
e^{-\int u^2} =   \int \langle \frac{\delta f}{\delta u},  
(u^2)_x\rangle 
e^{-\int u^2} =   -\int f\langle \frac{\delta }{\delta u}  
(u^2 )_x
e^{-\int u^2}\rangle  
\end{equation}
and
\begin{equation}
\langle \frac{\delta }{\delta u}  
(u^2 )_x
e^{-\int u^2}\rangle = 
\langle   
(2u_x - (u^2)_x 2u)\rangle 
e^{-\int u^2} .
\end{equation}
The last term vanishes because $(u^2)_x 2u = \frac23(u^3)_x$ and because of periodic boundary conditions any exact derivative integrates to zero: $\langle f_x\rangle
= \int_0^1 f =0$.  

Such an argument is known in physics \cite{S}.
Note that the problem is subtle, and requires 
an appropriate interpretation. In fact the result is not
correct for the
standard mathematical interpretation of the Burgers' flow as
the limit as $\epsilon\downarrow 0$ of $u^\epsilon_t
=2u^\epsilon u^\epsilon_x + \epsilon u^\epsilon_{xx}
$, as can be
checked with the Lax-Oleinik formula.  On the other
hand, the argument is rigorous for (\ref{sp}).

\section{Invariant measures for mK${\rm d}$V on $\mathbb T$}

Let $P_0$ denote Wiener measure on $\phi\in C(\mathbb T)$ conditioned to have
 $\int_{\mathbb T} \phi= 0$.   It can be derived from
 the standard circular Brownian motion $P$ on 
 $C(\mathbb T)$ defined as follows:  Condition a 
 standard Brownian motion $\beta (t)$, $t\in [0,1]$ 
 starting at $\beta(0)=x$ to have $\beta(1)=x$ as well, and now distribute $x$ on the real line according to Lebesgue measure.  $P_0$ is obtained from $P$ by conditioning on $\int_{\mathbb T} \phi=0$.

  Define 
$P^{(4)}_0$  to be the  measure absolutely continuous to $P_0$ given by
\begin{equation}
 P^{(4)}_0(B)
= Z^{-1}\int_{B} J(\phi)
e^{ -\frac12\int_{\mathbb T} \phi^4} 
dP_0,
\end{equation}
for Borel sets $B\subset C(\Bbb T)$ where $Z$ is the normalizing factor to make $ P^{(4)}_0$ a probability measure and
\begin{equation}
J(\phi) = (2\pi)^{-1/2}K(\phi)K(-\phi)e^{\frac12\left(\int\phi^2\right)^2}
\end{equation}
where
\begin{equation}
K(\phi)= \int_0^1  e^{2\Phi(x)} dx 
\end{equation}
and
\begin{equation}
\Phi(x)=\int_0^x \phi(y) dy.
\end{equation}

For smooth $g$ and $-\infty< t< \infty$, let $\phi(t)=\mathcal M_t g$ denote the (smooth) solution of the modified KdV (mKdV) equation,
\begin{equation}\label{mkdv}
\phi_t- 6\phi^2 \phi_x+\phi_{xxx} =0, \qquad \phi(0)=g.
\end{equation}

\begin{thm} (Kappeler and Topalov \cite{KT3})
$\mathcal M_t$ extends to a continuous group of nonlinear evolution operators \begin{equation}\bar{\mathcal M}_t: L^2(\Bbb T)\to L^2(\Bbb T).
\end{equation} \label{ktthm}
\end{thm}
Let\begin{equation}
  H(\phi) = \frac12 \int_{\Bbb T} \phi^4 + \phi_x^2. 
  \end{equation}
  mKdV can be written in Hamiltonian form,
  \begin{equation}
  \phi_t = \partial_x \frac{\delta H}{\delta\phi}.
  \end{equation}
$P_0^{(4)}$ gives rigorous meaning to the 
weighted Gibbs measure $J(\phi)e^{-H(\phi)}$ on ${\int\phi =0}$.  

\begin{thm} (Bourgain \cite{Bo}) $ P_0^{(4)}$ is invariant for
mKdV,
\begin{equation}
\bar{\mathcal M}_t^* P_0^{(4)}= P_0^{(4)}.
\end{equation}\label{theorem2.2} 
\end{thm}

\begin{proof} In fact what is proven in \cite{Bo} is that 
$Z^{-1}e^{-\frac12\int_{\mathbb T}\phi^4} dP$ is invariant for mKdV.  The main
obstacle at the time was a lack of well-posedness for mKdV on the support
$H^{1/2-}$ of the measure.  This statement follows with less work once one has the results  of Kappeler and Topalov proving
 well-posedness on a larger set (Theorem \ref{ktthm}).  
 
 We have in addition to show that $J(\phi)$ is a conserved quantity for
 mKdV.  It is well known that $\int_{\mathbb T}\phi^2$ is preserved.  So the 
 problem is reduced to showing that
 $
K(\phi)$ and $K(-\phi)$
are conserved.  Let $\phi(t)$ be a smooth solution of mKdV.  Note that
\begin{equation}
\partial_t\Phi = 2\phi^3 - \phi_{xx}.
\end{equation}
Hence
\begin{equation}
\partial_tK =  2\int_0^1 (2\phi^3 - \phi_{xx})  e^{2\Phi(x)}dx.  
\end{equation}
But integrating by parts we have, since $\phi$ is periodic and $\Phi_x=\phi$,
\begin{equation}
\int_0^1 \phi_{xx}  e^{2\Phi(x)}dx = 
-\int_0^1 2\phi_x \phi e^{2\Phi(x)}dx = 
-\int_0^1 (\phi^2)_x  e^{2\Phi(x)}dx
=\int_0^1 2\phi^3 e^{2\Phi(x)} dx.
\end{equation}
Therefore $\partial_tK(\phi(t)) =0$.  One can easily 
check with the analogous integration by parts
that $\partial_tK(-\phi(t))=0$.  

Now suppose $\bar{\mathcal M}_t\phi = \phi(t)$
with $\phi\in L^2(\mathbb T)$.  From Theorem \ref{ktthm} we have smooth $\phi_n$ with $\phi_n\to \phi$ and $\phi_n(t)\to \phi(t)$ in $L^2(\mathbb T)$.
\begin{equation}
K(\phi(t)) - K(\phi) 
= [K(\phi(t)) - K(\phi_n(t))] -[K(\phi_n)-K(\phi)]
\end{equation}
so if $K$ is a continuous
functions on $L^2(\mathbb T)$ then $K(\phi)$ and $K(-\phi)$ are conserved
by $\bar{\mathcal M}_t$.  To prove that
$K$ is continuous simply note that
\begin{equation}
|K(\phi)-K(\psi)| = |\int_0^1 e^{2\int_0^x \phi}
[ e^{ 2 \int_0^x \psi - \phi} -1] dx|\le
 e^{2\|\phi\|_{L^2(\mathbb T)}} 
[ e^{2\|\psi- \phi\|_{L^2(\mathbb T)}}-1]. 
\end{equation}
\end{proof}

\section{The Miura Transform on $L^2_0(\mathbb T)$}
\label{sec:miura}

The Miura transform $\phi\mapsto\phi_x + \phi^2 $ maps smooth solutions
of mKdV to smooth solutions of KdV.  It is basically two to one, and not onto.
But this is mostly a matter of the mean $\int_{\mathbb T} \phi$.  Since the mean
is conserved in both mKdV and KdV, it is more natural to consider the map corrected by
subtracting the mean.
The corrected Miura transform is
defined for smooth $\phi$ by,
\begin{equation}
\mu(\phi) = \phi_x + \phi^2 - \int_{\mathbb T} \phi^2.
\end{equation}

Let $L^2_0(\mathbb T)$ and $H^{-1}_0(\mathbb T)$ denote the subspaces
of $L^2(\mathbb T)$ and $H^{-1}(\mathbb T)$ with $\int_{\mathbb T}\phi=0$.

\begin{thm} (Kappeler and Topalov \cite{KT2})  The corrected Miura transform $\mu$
extends to a continuous map
\begin{equation}
\bar \mu : L^2_0\to H^{-1}_0
\end{equation}
which is one to one and onto.  $\bar\mu$
takes solutions $\phi$ of 
 mKdV
 (\ref{mkdv}) on $L^2(\Bbb T)$,
to solutions $u=\mu(\phi)$ of  KdV  (\ref{kdv})
on $H^{-1}(\Bbb T)$;
\begin{equation}\label{kt}
\bar{\mathcal S}_t\bar\mu = \bar\mu\bar{\mathcal M}_t.
\end{equation}
\label{theorem3.1}
  \end{thm}

{\it Remark.}  The Ricatti map is given by
\begin{equation}
r(\phi,\lambda)= \phi_x + \phi^2 +\lambda.
\end{equation}
Note that  Kappeler and Topalov use the term 
Ricatti map for $\mu= r(\phi, -\int_{\mathbb T} \phi^2)$.

\section{The Miura Transform on Wiener Space}
\label{sec:miura2}

\begin{thm} (Cambronero and McKean \cite{CM}) The corrected Miura transform $\bar \mu$ maps  $P_0^{(4)}$  into mean zero white
noise $Q_0$;
\begin{equation}\label{push}
\mu^{*}P_0^{(4)}=Q_0.
\end{equation}\label{theorem4.1}
\end{thm}

\begin{proof}  Let $\hat P_0^{(4)}$ be given by
\begin{equation}
\hat P^{(4)} (B)
= \frac1{\sqrt{2\pi}}\int_{(\phi,\lambda)\in B} K(\phi)K(-\phi)
e^{ -\frac12\int_{\mathbb T} (\phi^2+\lambda)^2  } 
dPd\lambda
\end{equation}
where $B$ is a Borel subset of $C(\mathbb T)\times \mathbb R$. Let 
$\hat r = (r, \int_{\mathbb T}\phi)$.  Let $\hat Q$ on $C(\mathbb T)\times \mathbb R$ be given by  $\hat Q= Q\times$ Lebesgue measure. What is actually proved in \cite{CM} is that
\begin{equation}
\hat r^*\hat P^{(4)} = \hat Q.
\end{equation}
(\ref{push}) is obtained by conditioning on $\lambda=-\int_{\mathbb T} \phi^2$ and 
$\int_{\mathbb T}\phi =0$.
\end{proof}

{\it Remark.}  There is a simple heuristic argument
explaining (\ref{push}).  Formally
\begin{equation}
dP_0^{(4)} = Z_1^{-1}K(\phi)K(-\phi) 
e^{-\frac12 \int_0^1 (\phi^2 + \phi' -\int_0^1\phi^2)^2}
dF(\phi),
\qquad dQ_0 = Z_2^{-1} e^{-\frac12\int_0^1 u^2}
dF(u)
\end{equation}
where $F$ is the (mythical) flat measure on $\int_0^1\phi=0$. Note that in the exponent of $dP_0^{(4)}$
we have assumed that integration by parts gives
$\int_0^1 \phi^2\phi' =0$.  Since the corrected
Miura transform $u = \phi^2 + \phi' -\int_0^1\phi^2$
the only mystery is the form of the Jacobian
$CK(\phi)K(-\phi)$.  Let $D$ be the map
$Df= f'$ and $\phi$ stand for the map of multiplication by $\phi$ with a subtraction to make
the result mean zero, $\phi f = \phi\cdot f - \int_0^1\phi\cdot f$.  The Jacobian is then,
\begin{equation}
f(\phi)= \det( 1+ 2\phi D^{-1})= \exp\{ {\rm Tr} \log (1+2\phi
D^{-1} )\}
\end{equation}
For fixed $x,y\in \T$ let $\partial_{xy}=
\frac{\partial}{\partial (\phi(y)-\phi(x))}$, i.e.~the G\^ ateaux derivative in the direction  $\delta_y-\delta_x$, $\partial_{xy}F(\phi) = 
\lim_{\epsilon\to 0} \epsilon^{-1}( F(\phi+ \epsilon(\delta_y-\delta_x)) - F(\phi)$. 
We have 
\begin{equation}
\partial_{xy} \log f(\phi)= \partial_{xy} {\rm Tr} \log (1+2\varphi
D^{-1}) = {\rm Tr} [ \{\partial_{xy}(1+2\phi
D^{-1})  \}\{ 1+2\phi
D^{-1})^{-1}\} ].
\end{equation}
If we let $G(x,y)$ denote the Green function 
of $D+ 2\phi$ this gives
$
\partial_{xy} \log f(\phi)= 2[G(y,y)-G(x,x)].
$.  It is not hard to compute the Green function 
with the result that
\begin{equation}
\partial_{xy} \log f(\phi)=\frac{ 2\int_x^y e^{-2\Phi}}{ \int_0^1 e^{-2\Phi}}-\frac{ 2\int_x^y e^{-2\Phi}}{ \int_0^1 e^{-2\Phi}}.
\end{equation}
The argument is completed by a straightforward verification that this is satisfied by $f(\phi)=K(\phi)K(-\phi)$.

The heuristic argument can be made rigorous by
taking finite dimensional approximations where this
set of equations actually identifies the determinant.  Since the
computations become exactly those of \cite{CM}, we do
not repeat them here.

\section{Proof of Theorem \ref{main}}
\label{sec:2proof}

\begin{equation} \bar{\mathcal S}_t^* Q_0  ~ {\stackrel{\rm Thm~ \ref{theorem4.1}}{=}}  ~\bar{\mathcal S}_t^*\mu^{*}P_0^{(4)}~ {\stackrel{\rm Thm~ \ref{theorem3.1}}{=}}~ \mu^* \bar{\mathcal M}_t^*P_0^{(4)}~{\stackrel{\rm Thm~ \ref{theorem2.2}}{=}}~ \mu^* P_0^{(4)}~{\stackrel{\rm Thm~ \ref{theorem4.1}}{=}}~ Q_0 
\end{equation}
%\begin{eqnarray*} \bar{\mathcal S}_t^* Q_0   &{\stackrel{\rm Thm~ \ref{theorem4.1}}{=}}~~   \bar{\mathcal S}_t^*\mu^{*}P_0^{(4)} &\qquad {\rm (Theorem ~ \ref{theorem4.1})}\\
%&=~~ \mu^* \bar{\mathcal M}_t^*P_0^{(4)} &\qquad {\rm (Theorem~ \ref{theorem3.1})}\\
%&=~~ \mu^* P_0^{(4)} &\qquad {\rm (Theorem~ \ref{theorem2.2})}\\
%&=~~ Q_0 &\qquad {\rm (Theorem~ \ref{theorem4.1})}
%\end{eqnarray*}
\medskip

\noindent{\sc
 Acknowledgements.}  
Thanks to K. Khanin, M. Goldstein and J. Colliander
for enlightening conversations.

\end{document}